\documentclass[a4paper,10pt]{article}

\newtheorem{thm}{Theorem}
\newtheorem{cor}[thm]{Corollary}

\newtheorem{lem}[thm]{Lemma}

\scrollmode
\def\qed{\hbox{\vrule height 7pt depth 0pt width 7pt}}
\def\cqfd{\hfill\penalty 500\kern 10pt\qed\medbreak}
 
\def \a{{\alpha}}
\def \b{{\beta}}

\def \O{{\Omega}}

\def \R{{\bf R}}

\def \E{{\bf E}}

\def \A{{\cal A}}

 \def \P{{\bf P}}

\def \qq{{\qquad}}
 \def \p{{\varphi}}
\def \noi{{\noindent}}

\def \e{{\varepsilon}}

\font\ph=cmcsc10  at  10 pt

\font\ph=cmcsc10  at  10 pt
\font\phh=cmcsc10  at  8 pt

 at 10,4 pt

 at 10,8 pt
 \def \lt {{\hbox{\vrule height 6pt depth 2pt width 0,7pt}\kern 1,1pt}}
\def \rt {{\kern 0,9pt\hbox{\vrule height 6pt depth 2pt width 0,7pt}}}

\rm
\def\ddate {\ifcase\month\or January\or
February\or March\or April\or May\or June\or July\or
August\or September\or October\or November\or December\fi\ {\the\day},
{\the\year}}
\font\phh=cmcsc10
at  8 pt
\scrollmode
\title{On Mean Values of Dirichlet Polynomials}
\author{
  Michel Weber}

\begin{document}

\maketitle

 \begin{abstract}  We show the following general lower bound valid for any positive integer $q$,  and arbitrary  reals $ \p_1,
\ldots,
\p_N  $ and non-negative reals   $a_1,\ldots, a_N$,    $$   c_{q}\Big(\sum_{n=1}^N  a_n^2\Big)^q  \le {1 \over 2T} \int_{|t |\le
T}   \Big|
\sum_{n=1}^N a_n e^{it\p_n}\Big|^{2q}dt. $$
 \end{abstract}


\section{Main Result}   The object of this short Note is to prove the following lower bound 
\begin{thm}  \label{t1} For any positive integer $q$, there exists a constant $c_q$, such that
for any reals $ \p_1,
\ldots,
\p_N  $,   any non-negative reals   $a_1,\ldots, a_N$, and any $T>0$,   
 $$   c_{q}\Big(\sum_{n=1}^N  a_n^2\Big)^q  \le {1 \over 2T} \int_{|t |\le
T}   \Big|
\sum_{n=1}^N a_n e^{it\p_n}\Big|^{2q}dt. $$
\end{thm}   
 The result is no longer true for arbitrary reals $a_1,\ldots, a_N$ as yields the case $\p_1=\ldots=\p_N$. It also follows that 
\begin{equation}  \label{e1} c  \Big(\sum_{n=1}^N  a_n^2\Big)^{1/2}  \le   \sup_{t\in \R}   \Big|
\sum_{n=1}^N a_n e^{it\p_n}\Big|  . 
\end{equation}
In the case $\p_n= \log n$, it is
known from
\cite{KQ} and \cite{Q3} that for any $(a_n)$
\begin{equation}\label{e2}
\sup_{t\in \R} \Big|
\sum_{n=0}^{N-1} a_n n^{it} \Big| \ge 
\alpha_1 { e^{\beta_1 \sqrt{\log N\log\log N} } \over
\sqrt{N}} 
 \Big(\sum_{n=0}^{N-1} |a_n|\Big)
\end{equation}
and for some
$(a_n)$
\begin{equation} \label{e3}
\sup_{t\in \R} \Big| \sum_{n=0}^{N-1} a_n n^{it} \Big| \le 
\alpha_2
{ e^{\beta_2 \sqrt{\log N\log\log N} } \over
\sqrt{N}} 
 \Big(\sum_{n=0}^{N-1} |a_n|\Big),
\end{equation}
with some universal constants
$\alpha_{1 },\alpha_{ 2},\beta_{1 },\beta_{ 2}$. 
  Then (\ref{e1}) is better than  (\ref{e2}) if for instance $a_n=  n^{-\a}$, $\a>1/2$, since   
$$ { e^{\beta_1 \sqrt{\log N\log\log N} } \over
\sqrt{N}} 
 \Big(\sum_{n=0}^{N-1} |a_n|\Big)\sim e^{\beta_1 \sqrt{\log N\log\log N} } )N^{{1\over 2}-\a}=o(1)\ll \Big(\sum_{n=1}^N 
a_n^2\Big)^{1/2}.
$$

\smallskip\par The $L^1$-case is related to well-known Ingham's inequality \cite{[In]}. We state the
sharper form   due to Mordell
\cite{[Mor]}: let
$0<\p_1<\ldots<\p_N$ and  let $\gamma$ be such that
$\displaystyle{\min_{1<n\le N}} \p_n-\p_{n-1}\ge \gamma>0$. Then  
\begin{equation}\sup_{n=1}^N|  a_n| \le {K\over T}\int_{-T}^T\Big|
\sum_{n=1}^N a_n e^{it\p_n}\Big| dt\qq {\rm with}\ T={\pi\over \gamma}, \label{4.1}\end{equation} 
where $K\le 1$.

\smallskip Further  with no restriction, one
always have
\begin{equation}\sup_{n=1}^N |a_n|\le \limsup_{T\to \infty}{1\over 2T}\int_{-T}^T\Big|
\sum_{n=1}^N a_n e^{it\p_n}\Big| dt \le \sup_{t\in \R}\Big|
\sum_{n=1}^N a_n e^{it\p_n}\Big|, \label{4.2}\end{equation}
a very familiar inequality    in the theory of uniformly almost periodic functions. 
  See also   \cite{[B]} where the more complicated inequality is established:
\begin{equation}|a_n| \le {1\over \prod_{j=0}^{n-1}\cos({\pi\p_j\over 2\p_n})\cdot\prod_{n+1}^N\cos({\pi\p_n\over 2\p_j})}\cdot\sup_{|t|\le {\pi\over
2}({n\over \p_n}+\sum_{j=n+1}^N{1\over \p_j})}\Big|
\sum_{n=1}^N a_n e^{it\p_n}\Big|. 
\end{equation}    In   particular, if   $ \p_1,
\ldots,
\p_N $  are   linearly independent, and  $T$ is  large enough, then  \begin{equation}  b_{q}\Big(\sum_{n=1}^N  a_n^2\Big)^q  \le { 1\over 2T} \int_{|t |\le
T}   \Big|
\sum_{n=1}^N a_n e^{it\p_n}\Big|^{2q}dt \le B_{q}\Big(\sum_{n=1}^N  a_n^2\Big)^q,  \label{4.3}\end{equation} 
   holds    for any
nonnegative reals 
$a_1,\ldots, a_N$      
and   $b_q$, $B_q$    depend  on $q$ only. 
  \bigskip\par
      The proof of Theorem \ref{t1}  relies upon the following  lemma, which just
generalizes a useful majorization argument (\cite{[Mon]}, p.131) to arbitrary even powers. 
\begin{lem} \label{lm}Let $q$ be
any positive integer. Let     
$c_1,\ldots, c_N$ be   complex numbers and nonnegative reals $a_1,\ldots, a_N$  such that $|c_n|\le a_n$, 
$n=1,\ldots, N$. Then for any reals  $ T,T_0$ with
$T>0$ 
$$ \int_{|t- T_0|\le T}\Big| \sum_{n=1}^N c_ne^{it\p_n}\Big|^{2q}dt \le 3 \int_{|t |\le T}   \Big|
\sum_{n=1}^N a_ne^{it\p_n}\Big|^{2q}dt. $$
\end{lem}
 {\it Proof.}    Let 
$$ K_T(t)=K_T(|t|)=\big( 1-|t|/T)\chi_{\{|t|\le T\}} $$
Observe that for any reals $t,H$
\begin{eqnarray*}a)\quad & &K_T(t-H) =\big( 1-|t-H|/T)\chi_{\{|t-H|\le T\}}
 \cr b)\quad&&\chi_{\{|t-H|\le T\}} \le K_T(t-H)+K_T(t-H+T)+ K_T(t-H-T)\cr  
c)\quad &&\widehat{ K}_T(u)  ={1\over T}\big({\sin  Tu\over   u  }\big)^2 \ge 0,\qq\hbox{for all real $u$}  .
\end{eqnarray*}
  Suppose that $|c_n|\le a_n$ for $n=1,\ldots, N$. From 
\begin{eqnarray}\Big( \sum_{n=1}^N c_ne^{it\p_n}\Big)^{ q} &=&
\sum_{k_1+\ldots+ k_N= q} \Big({q!\over k_1!\ldots k_N!}\Big)\prod_{n=1}^Nc_n^{ k_n} e^{i t k_n\p_n }  
 .
\end{eqnarray} 
 and $$\Big| \sum_{n=1}^N c_ne^{it\p_n}\Big|^{2 q} = \sum_{k_1+\ldots+ k_N= q\atop h_1+\ldots+ h_N= q} \Big({(q!)^2\over k_1!h_1!\ldots k_N!h_N!} 
\Big)\prod_{n=1}^Nc_n^{ k_n}{\overline{c }_n} ^{ h_n} e^{i t (k_n-h_n)\p_n }    $$  
we get 
 \begin{eqnarray*}&&\int_\R K_T (t-H) \Big| \sum_{n=1}^N c_ne^{it\p_n}\Big|^{2q}dt \cr 
&=&\sum_{k_1+\ldots+ k_N= q\atop h_1+\ldots+ h_N= q} 
  {(q!)^2\over k_1!h_1!\ldots k_N!h_N!} 
\ \prod_{n=1}^Nc_n^{ k_n}{\overline{c }_n} ^{ h_n}    \int_\R
K_T(t-H)e^{it\sum_{n=1}^N (k_n-h_n)\p_n}dt\cr  
&=&\sum_{k_1+\ldots+ k_N= q\atop h_1+\ldots+ h_N= q} 
  {(q!)^2\over k_1!h_1!\ldots k_N!h_N!} 
\ \prod_{n=1}^Nc_n^{ k_n}{\overline{c }_n} ^{ h_n}    \int_\R
K_T(s)e^{i(s+H)\sum_{n=1}^N (k_n-h_n)\p_n}ds 
\end{eqnarray*}
$$\quad =\sum_{k_1+\ldots+ k_N= q\atop h_1+\ldots+ h_N= q} 
 {(q!)^2\over k_1!h_1!\ldots k_N!h_N!} 
  \prod_{n=1}^N(c_ne^{i H     \p_n})^{ k_n} \overline{(c_ne^{i H     \p_n}})^{ h_n}  \widehat{
K}_T\Big(\sum_{n=1}^N  (k_n-h_n)\p_n \Big)
$$
   \begin{eqnarray*}&\le & 
\sum_{k_1+\ldots+ k_N= q\atop h_1+\ldots+ h_N= q} 
 {(q!)^2\over k_1!h_1!\ldots k_N!h_N!} 
\ \prod_{n=1}^N a_n ^{ k_n+ h_n} 
  \widehat{
K}_T\Big(\sum_{n=1}^N  (k_n-h_n)\p_n \Big)\cr &
= &
  \int_\R  {
K}_T(t) \bigg[\sum_{k_1+\ldots+ k_N= q\atop h_1+\ldots+ h_N= q} 
 {(q!)^2\over k_1!h_1!\ldots k_N!h_N!} 
  \prod_{n=1}^N a_n ^{ k_n+ h_n}   e^{it\sum_{n=1}^N  (k_n-h_n)\p_n }  \bigg]  dt\cr &=&\int_\R K_T(t )\Big|
\sum_{n=1}^N a_ne^{it\p_n}\Big|^{2q} dt. 
\end{eqnarray*}
Hence, if $|c_n|\le a_n$ for $n=1,\ldots, N$  
\begin{equation}   \int_\R K_T(t- H)\Big| \sum_{n=1}^N
c_ne^{it\p_n}\Big|^{2q}dt\le \int_\R K_T(t )\Big|
\sum_{n=1}^N a_ne^{it\p_n}\Big|^{2q}dt.\label{4.5}\end{equation}

By applying Lemma \ref{lm} with $H=0, T_0, -T_0$,  and using b), we get 
\begin{eqnarray}&&\int_{|t- T_0|\le T}  \Big| \sum_{n=1}^N c_ne^{it\p_n}\Big|^{2q}dt\cr &\le &\int_\R \Big( K_T(t-T_0)+K_T(t-T_0+T)+
K_T(t-T_0-T)\Big)\Big|
\sum_{n=1}^N c_ne^{it\p_n}\Big|^{2q}dt \cr & \le& 3\int_\R   K_T(t ) \Big|
\sum_{n=1}^N a_ne^{it\p_n}\Big|^{2q}dt\le 3\int_{|t |\le T}    \Big|
\sum_{n=1}^N a_ne^{it\p_n}\Big|^{2q}dt. 
\end{eqnarray}
\cqfd
  
The  proof of Theorem \ref{t1} is now achieved as follows.   First recall the   Khintchin-Kahane inequalities \cite{KS}.  Let $\{
\e_i, 1\le i\le N\}$ be   independent Rademacher random variables,  thus satisfying $\P\{ \e_i =\pm 1\} =1/2$, if $(\O, \A,
\P)$  denotes the underlying  basic probability.    Then for any $0<p<\infty$, there
exist positive finite constants
$c_p$,
$C_p$ depending on
$p$ only, such that for any   sequence $\{ a_i, 1\le i\le N\}$ of real numbers
\begin{equation} \label{Khintchin} c_p \Big(\sum_{i=1}^N a_i^2\Big)^{1/2}\le \Big\|\sum_{i=1}^N a_i\e_i\Big\|_p\le C_p \Big(\sum_{i=1}^N
a_i^2\Big)^{1/2}.
\end{equation}  This remains true for complex $a_n$. If $a_n =\a_n +i\b_n$, then 
 \begin{eqnarray*}\Big\|\sum_{j=1}^N a_j\e_j\Big\|_p^p&=&\E \Big|\sum_{j=1}^N \a_j\e_j+i\sum_{j=1}^N \b_j\e_j\Big|^p=\E\Big(
\Big|\sum_{j=1}^N
\a_j\e_j\Big|^2+  \Big|\sum_{j=1}^N \b_j\e_j\Big|^2\Big)^{p/2}\cr &\le& 2^{(p/2)-1}\Big(\E \Big|\sum_{j=1}^N
\a_j\e_j \Big|^p+\E \Big| \sum_{j=1}^N \b_j\e_j\Big|^p\Big), \end{eqnarray*}
where we have denoted by $\E$   the  corresponding expectation symbol. Thus, since $\sqrt A+\sqrt B\le \sqrt{2(A+B)}$, $A,B\ge 0$,
\begin{eqnarray*}\Big\|\sum_{j=1}^N a_j\e_j\Big\|_p  &\le& 2^{1/2-1/p}C_p\Big[\big(\sum_{j=1}^N \a_j^2\big)^{1/2} +\big(\sum_{j=1}^N
\b_j^2\big)^{1/2}\Big]
\cr &\le &2^{1-1/p}C_p \Big(\sum_{j=1}^N (\a_j^2 +   \b_j^2)\Big)^{1/2}   =   C'_p \Big(\sum_{j=1}^N |a_j|^2 
\Big)^{1/2}.\end{eqnarray*} 
Conversely, from
\begin{eqnarray*}\Big\|\sum_{j=1}^N a_j\e_j\Big\|_p^p&=&\E\Big( \Big|\sum_{j=1}^N
\a_j\e_j\Big|^2+  \Big|\sum_{j=1}^N \b_j\e_j\Big|^2\Big)^{p/2}\cr &\ge&  \max\Big(\E  \Big|\sum_{j=1}^N
\a_j\e_j\Big| ^{p }, \E   \Big|\sum_{j=1}^N \b_j\e_j\Big| ^{p }\Big),\end{eqnarray*}
we get 
\begin{eqnarray*}\Big\|\sum_{j=1}^N a_j\e_j\Big\|_p
&\ge&\max\Big(\big\|\sum_{j=1}^N \a_j\e_j\big\|_p, \big\|\sum_{j=1}^N \b_j\e_j\big\|_p\Big)\cr &\ge  &c_p\max\Big(\big(\sum_{j=1}^N
|\a_j|^2  \big)^{1/2},
\big(\sum_{j=1}^N |\b_j|^2  \big)^{1/2}\Big) 
\cr &\ge &{c_p\over 2}   \Big( \sum_{j=1}^N
(|\a_j|^2 + |\b_j|^2   ) \Big)^{1/2}=c'_p \Big( \sum_{j=1}^N
 | a_j|^2  \Big)^{1/2}. \end{eqnarray*}
Now choose $c_n=\e_na_n
$. Taking expectation in inequality of Lemma 2.1, and using Fubini's Theorem, gives  
\begin{equation} \int_{|t |\le T}  \E\,  \Big|
\sum_{n=1}^N \e_na_ne^{i t \p_n} \Big|^{2q} \le 3  \int_{|t |\le
T}   \Big|
\sum_{n=1}^N a_n e^{it\p_n}\Big|^{2q}dt. \label{4.7}\end{equation}
By (\ref{Khintchin}) we have
\begin{equation}c_{q}\Big(\sum_{n=1}^N  a_n^2\Big)^q\le \E\, \Big| \sum_{n=1}^N \e_na_ne^{it\p_n}\Big|^{2q}\le
 C_q\Big(\sum_{n=1}^N  a_n^2\Big)^q.\label{4.10}\end{equation}
By reporting 
 \begin{equation} 2T c_{q}\Big(\sum_{n=1}^N  a_n^2\Big)^q  \le 3  \int_{|t |\le
T}   \Big|
\sum_{n=1}^N a_n e^{it\p_n}\Big|^{2q}dt,\label{4.11}\end{equation}which proves our claim. 

 \cqfd
\section{Application}

  We shall deduce from Theorem \ref{t1} the following lower bound.
\begin{cor} For every $N$, $T$ and $\nu$ $$c_\nu\log^{\nu^2}  N\le
{1
\over 2T}
\int_{| t |\le T}   \Big|
\sum_{n=1}^N {1\over n^{{1\over 2}+i t}} \Big|^{2\nu}d t. $$
\end{cor}
  In relation with this is Ramachandra's  well-known lower bound   (see \cite{Iv} section 9.5, to which we
also refer for the estimates used in the proof)    
\begin{equation} \label{ram}   c_{\nu}(\log   T)^{\nu^2} \le {1 \over 2T} \int_{|t |\le T}   \big|
\zeta({1\over 2}+it ) \big|^{2\nu}dt. 
\end{equation} 
{\it Proof.} Apply Theorem \ref{t1} with $q=2$ to  the sum
$$
 \Big(\sum_{n=1}^N  {1\over
n^{ {1\over 2}+it  } }\Big)^{\nu}:=\sum_{m=1}^{N^\nu}{b_m\over m^{ {1\over 2}+it  } }, 
$$
where 
$$b_m=\#\big\{(n_j)_{j\le \nu}\, ;\, n_j\le
N:m=\prod_{j\le \nu} n_j\big\}   .  $$ Thus for all $N$ and $T$
 $$c_\nu\ \sum_{m=1}^{N^\nu}{b^2_m\over m  } \le
{1
\over 2T}
\int_{| t |\le T}   \Big|
\sum_{n=1}^N {1\over n^{{1\over 2}+i t}} \Big|^{2\nu}d t. $$
But if $m\le N$, $b_m=d_{ \nu}(m)$ where $d_\nu(m)$ denotes the number of representations of   $m$ as a product of
$\nu$ factors, and we know that
$$\sum_{m\le x} {d_{ \nu}^2(m)\over m}= (C_\nu +o(1))\log^{\nu^2} x. $$
Thus 
$$\sum_{m=1}^{N^\nu}{b^2_m\over m  }\ge \sum_{m=1}^{N }{b^2_m\over m  } \ge c_\nu \log^{\nu^2} N$$
Henceforth $$c_\nu\log^{\nu^2}  N\le
{1
\over 2T}
\int_{| t |\le T}   \Big|
\sum_{n=1}^N {1\over n^{{1\over 2}+i t}} \Big|^{2\nu}d t. $$
\cqfd 
 
\noi {\it Acknowlegments.}    I thank Professor Aleksandar Ivi\'c for useful remarks.

 {
}
\bigskip\par
 \noi {\phh Michel  Weber, \noi  Math\'ematique (IRMA),
Universit\'e Louis-Pasteur et C.N.R.S.,   7  rue Ren\'e Descartes,
67084 Strasbourg Cedex, France.
\par\noindent
E-mail: \  \tt weber@math.u-strasbg.fr} 

\baselineskip 12pt
\begin{thebibliography}{99}
\bibitem{[B]} {Binmore K.G.} [1966]: {\sl A trigonometric inequality}, J. London Math. Soc. {\bf 41}, 693--696.
     \bibitem{[In]}  {Ingham A.E.}  [1950]: {\sl A further note on trigonometrical inequalities}, Proc. Cambridge Philos. Soc. {\bf
46}, 535-537. 
\bibitem{Iv} {Ivi\'c A.} [1985] {\sl The Riemann Zeta-function},
Wiley-Interscience Publication, J. Wiley\&Sons, New-York.
\bibitem{KS} {Kashin B.S., Saakyan A.A.} [1989]    {\sl Orthogonal Series}, Translations of
Mathematical Monographs {\bf 75}, American Math. Soc. 

\bibitem{KQ} {Konyagin S.V., Queff\'elec H.} [2001/2002] {\sl
 The translation ${1\over 2}$ in the theory of Dirichlet series},
 Real Anal. Exchange {\bf 27}(1), 155--176.

 
 \bibitem{[Mon]} {\ph Montgomery H.} [1993]:  {\sl  Ten lectures on the interface between analytic number theory and harmonic analysis},
Conference Board of the Math. Sciences, Regional Conference Series in Math. {\bf 84}.
 \bibitem{[Mor]} {\ph Mordell I.J.} [1957]: {\sl On Ingham's trigonometric inequality}, Illinois J. Math. {\bf 1},
214--216. 
\bibitem{Q3} {Queff\'elec H.} [1995] {\sl H. Bohr's vision of
 ordinary Dirichlet series; old and new results},  J. Analysis  {\bf 3},
 p.43-60.
\bibitem{R} {Ramachandra K.} [1995] {\sl On the Mean-Value and Omega-Theorems for the Riemann Zeta-Function},  Tata Institute of
Fundamental Research, Bombay, Springer Verlag  Berlin, Heidelberg, New-York, Tokyo, vii+167p.  \end{thebibliography}
\end{document}